\input amstex
 \documentstyle{amsppt}
\NoBlackBoxes
\magnification1200
\pagewidth{6.5 true in} 
\pageheight{9.25 true in}
\topmatter
\title
The distribution of prime numbers 
\endtitle 
\author K. Soundararajan 
\endauthor 
\thanks 
The author is partially supported by the National Science Foundation.
\endthanks
\address 
Department of Mathematics, University of Michigan, Ann Arbor, MI 48109, USA
\endaddress
\email 
ksound\@ \!\!umich.edu
\endemail
\endtopmatter
\def\lam{\lambda}
\def\Yildirim{Y{\i}ld{\i}r{\i}m\ }
\document
\loadbold
 \document

 What follows is an expanded version of my lectures at the NATO School on Equidistribution.  
 I have tried to keep the informal style of the lectures.  In particular, I have sometimes oversimplified 
 matters in order to convey the spirit of an argument.   

 \head Lecture 1:  The Cram{\' e}r model and gaps between consecutive primes \endhead

 \noindent The prime number theorem tells us that $\pi(x)$, 
 the number of primes below $x$, is $\sim x/\log x$.  Equivalently, if $p_n$ 
 denotes the $n$-th smallest prime number then $p_n \sim n \log n$.  
 What is the distribution of the gaps between 
 consecutive primes, $p_{n+1}-p_n$?  

We have just seen that 
$p_{n+1}-p_n$ is approximately $\log n$ ``on average".  How often do we get a gap of size $2\log n$, say; or of 
size $\frac 12 \log n$?  One way to make this question precise is to fix an 
interval $[\alpha,\beta]$ (with $0\le \alpha<\beta$) and ask for 
$$
\lim_{N\to \infty} \frac 1N \#\Big\{ 2\le n\le N: \  \ \frac{p_{n+1}-p_n}{\log n} 
\in [\alpha,\beta]\Big\}. \tag{1.1}
$$
Does this limit exist, and if so what does it equal?

Here is another way to formulate this question.  Consider intervals 
of the form $[n,n+\log n]$ as $n$ ranges over integers up to $N$.  
On average such an interval contains one prime.  But of course 
some intervals may not contain any prime, and others may contain several.  
Given a non-negative integer $k$, how often does such an interval 
contain exactly $k$ primes?  What is 
$$
\lim_{N\to \infty} \frac{1}{N} \#\{ n\le N: \ \ \pi(n+\log n)-\pi(n) = k\}? \tag{1.2}
$$
Or more generally, for a fixed real number $\lambda >0$ we may ask for 
$$
\lim_{N\to \infty} \frac{1}{N} \#\{ n\le N: \ \ \pi(n+\lambda\log n)-\pi(n) = k\}? \tag{1.3}
$$

 In this lecture we will describe the conjectured answers to these questions, but we 
 confess at the outset that no one knows how to prove those conjectures.  
 While conjecturing the prime number theorem, Gauss stated that 
 the `density of primes' around $x$ should be $1/\log x$.   He based his conjecture on
 extensive numerical investigations.  In particular he divides the numbers 
 up to three million into intervals of length $100$ (a ``centad") and meticulously 
 tabulates the number of centads with no primes, exactly one prime etc.\footnote{We refer the 
 reader to www.math.princeton.edu/\~ \,\!ytschink/.gauss for scans of Gauss's manuscripts showing 
 these calculations.}  While he does not 
 seem to make a synthesis of his results (except to conjecture the prime number theorem) 
 it seems clear that he was seeking to understand a question like (1.3).    It was left 
 to Harald Cram{\' e}r [5] to set Gauss's work on a probabilistic footing.   

 \proclaim{Cram{\' e}r's model}  The primes behave like independent random variables 
 $X(n)$ ($n\ge 3$) 
with $X(n) = 1$ (the number $n$ is `prime') with probability $1/\log n$, 
and $X(n)=0$ (the number $n$ is `composite') with probability $1-1/\log n$.  
\endproclaim

 Let us suppose that the primes behave like a typical sequence in this random model, and 
 answer questions (1.1) and (1.3).   We want  
 the `probability' that $p_{n+1}-p_n$ lies between $\alpha \log n$ 
and $\beta \log n$.  Thus, given the prime $p_n$, we want $p_n +1$, 
$\ldots$, $p_{n}+h-1$ to be composite, and $p_{n}+h$ to be prime, 
where $\alpha \log n \le h\le \beta \log n$.  According to Cram{\' e}r's model, 
this occurs with 
probability 
$$
\align
\sum_{\alpha  \log n\le h \le \beta \log n} 
\prod_{j=1}^{h-1} \Big( 1- \frac{1}{\log (p_n + j)}\Big) \frac{1}{\log (p_{n}+h)} 
&\sim \sum_{\alpha  \log n\le h \le \beta \log n} 
\Big( 1-\frac{1}{\log n}\Big)^{h-1} \frac{1}{\log n} 
\\
\endalign
$$
since $\log (p_n+j) \sim \log n$ as $p_n\sim n\log n$ and $j \ll \log n$.   
This is 
$$
\sim \sum_{\alpha  \le h/\log n \le \beta} 
e^{-h/\log n} \frac{1}{\log n} \sim \int_{\alpha}^{\beta} e^{-t} dt,
$$
for large $n$, because the LHS looks like a Riemann sum approximation 
to the integral in the RHS. 
This is the conjectured answer to question (1.1): the probability ``density'' of finding $p_{n+1}-p_n$ 
close to $t\log n$ is $e^{-t}$.   This is an example of what is known as a ``Poisson process" 
in the probability literature, see Feller [7].  

\proclaim{Exercise 1}  Show similarly that the Cram{\' e}r model predicts that the answer to question (1.3) 
is $\frac{\lam^{k} }{k!} e^{-\lam}$.  This is the Poisson distribution with parameter $\lambda$.  
\endproclaim

The reader may well object that these predictions are dubious:  clearly the probability that 
$n$ and $n+1$ are both primes must be zero, but the Cram{\' e}r model assigns this event 
a probability $1/(\log n \log (n+1))$.    More generally, suppose we are given a set ${\Cal H}= 
\{h_1, \ldots, h_k\}$ of $k$ distinct integers, and we ask for the 
number of integers $n\le x$ with $n+h_1$, $n+h_2$, $\ldots$, $n+h_k$ all 
being prime.  The Cram{\' e}r model would predict an answer of $\sim x/(\log x)^k$, but 
clearly we must take into account arithmetic properties of the set ${\Cal H}$.  For 
example, if there were a prime $p$ such that the integers $h_1$, $\ldots$, $h_k$ 
occupied all the residue classes $\pmod p$ then the integers $n+h_1$, $\ldots$, 
$n+h_k$ would also occupy all the residue classes $\pmod p$.   In particular 
one of these numbers would be a multiple of $p$, and so there can only be finitely 
many values of $n$ with $n+h_1$, $\ldots$, $n+h_k$ all being prime.  

In [17] Hardy and Littlewood proposed the prime $k$-tuple conjecture 
that 
$$
\#\{n \le x:\  n+h_1, n+h_2, \ldots, n+h_k \text{ prime} \} 
\sim {\frak S}({\Cal H}) \frac{x}{(\log x)^k}, \tag{1.4} 
$$
for a certain constant ${\frak S}({\Cal H})$ called the `singular series.'  The constant 
${\frak S}({\Cal H})$ equals $0$ if the elements $h_1$, $\ldots$, $h_k$ occupy 
a complete set of residue classes $\pmod p$ for some prime $p$, and ${\frak S}({\Cal H})$ 
is positive otherwise.  We will describe this conjecture in more detail below.  
The aim of this lecture is to describe a beautiful calculation of
Gallagher [9] which shows that the Hardy-Littlewood conjecture (1.4) 
implies the same distribution of gaps between primes predicted by the Cram{\' e}r 
random model.  The crux of his proof is that although ${\frak S}({\Cal H})$ 
is not always $1$ (as the Cram{\' e}r model would have), it is $\sim 1$ on average over 
all $k$-element sets ${\Cal H}$ with the $h_j \le h$.

\noindent{\bf The Hardy-Littlewood Conjecture.}  We now motivate the Hardy-Littlewood 
conjecture (1.4) and describe the singular series ${\frak S}({\Cal H})$ that arises there.   As a toy 
model for prime numbers let us fix an integer $q$ and consider the reduced 
residue classes $\pmod q$.  Out of the $q$ 
total residue classes, there are $\phi(q)$ reduced classes, and we may think of 
$\phi(q)/q$ as the `probability' of a class being reduced.  Now suppose we 
are given the set ${\Cal H}=\{h_1, \ldots, h_k\}$ and we ask for the number of 
$n\pmod q$ such that $n+h_1$, $\ldots$, $n+h_k$ are all coprime to $q$. 
For convenience, let us just think of square-free $q$.  
If these $k$ events were independent then the answer would be $q (\phi(q)/q)^k$.  
The correct answer is a little different:  for each prime $p$ that divides $q$ 
we need $n$ to avoid the residue classes $-h_1$, $\ldots$, $-h_k \pmod p$.  
Let $\nu_{\Cal H}(p)$ denote the number of distinct residue classes 
occupied by ${\Cal H} \pmod p$.   Thus $n$ must lie in one of 
$p-\nu_{{\Cal H}}(p)$ residue classes $\pmod p$.  Using the chinese remainder theorem 
we see easily that the correct answer is
$$
\prod_{p|q} (p- \nu_{\Cal H}(p)) 
= q \prod_{p|q} \Big( 1-\frac{\nu_{\Cal H}(p)}{p}\Big) 
= q \Big(\frac{\phi(q)}{q}\Big)^k \prod_{p|q} 
\Big( 1- \frac{\nu_{\Cal H}(p)}{p} \Big) \Big(1-\frac{1}{p}\Big)^{-k}.
$$

Let us write ${\frak S}({\Cal H};q)= \prod_{p|q} (1-\frac{\nu_{\Cal H}(p)}{p})(1-\frac 1p)^{-k}$.  
We have seen that the answer for the number of $n \pmod q$ 
with $n+h_1$, $\ldots$, $n+h_k$ all being coprime to $q$ involves correcting 
the guess $q (\phi(q)/q))^k$ by the factor ${\frak S}({\Cal H};q)$ which keeps track 
of the arithmetic properties of the set ${\Cal H}$.  Now let us consider what 
happens when we take $q =q_{\ell} =\prod_{p\le \ell} p$ and let $\ell$ go to infinity.  
As $\ell \to \infty$ we see that 
$$
{\frak S}({\Cal H};q_{\ell}) \to \prod_{p} \Big(1-\frac{\nu_{\Cal H}(p)}{p}\Big) \Big(1-
\frac{1}{p}\Big)^{-k}.
$$
The infinite product above converges because if $p$ is larger than all the $h_j$'s 
then $\nu_{\Cal H}(p)=k$ and so $(1-\nu_{\Cal H}(p)/p) (1-1/p)^{-k} = 
(1-k/p) (1-1/p)^{-k}=1+O(p^{-2})$.  This infinite product is the singular series\footnote{The 
term arises from Hardy and Littlewood's original derivation of their conjecture using 
the circle method.  Here ${\frak S}({\Cal H})$ arose as a series rather than the 
product given above.}:
$$
{\frak S}({\Cal H}) := \prod_{p} \Big(1- \frac{\nu_{\Cal H}(p)}{p}\Big) \Big(1-\frac 1p\Big)^{-k}. 
\tag{1.5}
$$
Further if $n+h_1$, $\ldots$, $n+h_k$ are coprime to $q_{\ell}$ with $\ell$ large, then 
they have no small prime divisors, and may reasonably be viewed as a kind of 
approximation to primes.  Thus in formulating a conjecture on the number of $n\le x$ 
with $n+h_1$, $\ldots$, $n+h_k$ being prime, a natural guess is to take the 
random answer $x/(\log x)^k$, and multiply it by the arithmetical correction factor ${\frak S}({\Cal H})$. 
This is precisely the Hardy-Littlewood conjecture (1.4).  
It is immediate from (1.5)  that ${\frak S}({\Cal H})=0$ if and only if ${\Cal H}$ 
exhausts a compete set of residue classes $\pmod p$ for some $p$.

\noindent{\bf Gallagher's calculation.}  We will now describe Gallagher's argument, using 
the Hardy-Littlewood conjecture (1.4) to justify the distribution of gaps between primes 
predicted by the Cram{\' e}r model.   The precise problem we consider is a close variant 
of question (1.3).  Let $\lambda$ be a positive real number, and let $N$ be large.  We set 
$h=\lambda \log N$ and seek to understand the distribution of $\pi(n+h)-\pi(n)$ 
as $n$ varies over the natural numbers below $N$.   To understand this quantity, we 
consider the moments
$$
\frac{1}{N} \sum_{n\le N} (\pi(n+h)-\pi(n))^r = \frac 1N \sum_{n\le N} \Big( \sum\Sb \ell=1 \\ 
n+\ell \text{ prime} \endSb^h 1\Big)^r,  \tag{1.6}
$$
where $r$ is a natural number.   If the Cram{\' e}r prediction is right, then we may expect these 
moments to be approximately 
$$
\frac{1}{N} {\Bbb E}\Big(\sum_{2\le n\le N} \Big( \sum_{\ell=1}^{h} X(n+\ell) \Big)^r\Big),\tag{1.7}
$$
where ${\Bbb E}$ denotes expectation, and the $X(n)$'s are independent random variables 
as in Cram{\' e}r's model.  If these moments are roughly equal for $r\le R$ for any $R=R(N)$ 
tending to infinity, then we would know that 
$\pi(n+h)-\pi(n)$ has a Poisson distribution with parameter $\lambda$.   This is 
because of a well-known principle from probability, that nice distributions 
including the Poisson distribution are determined by their moments.

Let us expand out the $r$-th powers in (1.6) and (1.7).  We then 
get numbers $\ell_1$, $\ldots$, $\ell_r$ below $h$ not necessarily distinct 
and would like to understand how often $n+\ell_1$, $\ldots$, $n+\ell_r$ 
are all prime (for (1.6)), or to understand ${\Bbb E}(X(n+\ell_1)\cdots X(n+\ell_r))$ 
(for (1.7)).   Let us suppose that there are exactly $k$ distinct numbers among 
the $\ell_1$, $\ldots$, $\ell_r$ and write these distinct numbers as 
$(1\le )h_1 < h_2 <\ldots < h_k (\le h)$.   The number of choices for $\ell_1$, 
$\ldots$, $\ell_r$ that lead to the same ordered set of distinct numbers 
$h_1$, $\ldots$, $h_k$ is the number of different ways of mapping 
$\{1,2,\ldots, r\}$ onto $\{1,\ldots, k\}$; let us denote this\footnote{This is 
a `Stirling number of the second kind.'} 
 by $\sigma(r,k)$.  Thus we see that (1.6) may be written as 
 $$
 \sum_{k=1}^{r} \sigma(r,k) \sum\Sb 1\le  h_1 < h_2 <\ldots<h_k \le h\endSb 
 \ \ \Big(\frac{1}{N} \sum\Sb n\le N \\ n+h_1, \ldots, n+h_k \text{ prime}\endSb 1\Big), \tag{1.8}
 $$
 while (1.7) may be written as 
 $$
  \sum_{k=1}^{r} \sigma(r,k) \sum\Sb 1\le  h_1 < h_2 <\ldots<h_k \le h\endSb 
\ \ \Big(\frac{1}{N} \sum_{2\le n\le N} {\Bbb E}(X(n+h_1)\cdots X(n+h_k))\Big). \tag{1.9}
$$
Since   
the same quantity appears in both (1.8) and (1.9) and is non-negative, we don't need to worry about what 
$\sigma(r,k)$ is.

Invoking the Hardy-Littlewood conjecture (1.4)\footnote{Precisely, we need this conjecture 
uniformly for all  $h_1$, $\ldots$, $h_k$ below $h$, and for all $k\le R$ with 
$R=R(N)$ tending slowly to infinity.}  we get that (1.8) is 
$$
\sim \sum_{k=1}^{r} \frac{\sigma(r,k)}{(\log N)^k} 
\sum\Sb 1\le  h_1 < h_2 <\ldots<h_k \le h\endSb {\frak S}(\{h_1,\ldots,h_k\}).
$$
Clearly the quantity in (1.9) is 
$$
\sim  \sum_{k=1}^{r} \frac{\sigma(r,k)}{(\log N)^k} 
\sum\Sb 1\le  h_1 < h_2 <\ldots<h_k \le h\endSb 1.
$$
Thus, to show that (1.6) and (1.7) are approximately equal, we need only 
show that 
$$
\sum\Sb 1\le  h_1 < h_2 <\ldots<h_k \le h\endSb {\frak S}(\{h_1,\ldots,h_k\}) 
\sim \sum\Sb 1\le  h_1 < h_2 <\ldots<h_k \le h\endSb 1. \tag{1.10}
$$
This is Gallagher's crucial result in [9].  It shows that although the Hardy-Littlewood 
probabilities are different from the Cram{\' e}r probabilities, on average they are 
roughly equal.  This explains why the Cram{\' e}r model makes accurate 
predictions for the distribution of primes in such short intervals.  

\proclaim {Exercise 2}  For a prime $p$ put ${\frak S}({\Cal H};p) = (1-\nu_{\Cal H}(p)/p)(1-1/p)^{-k}$. 
Prove that as $h\to \infty$ 
$$
\sum\Sb  1\le  h_1 < h_2 <\ldots<h_k \le h\endSb {\frak S}({\Cal H};p) 
\sim  \sum\Sb 1\le  h_1 < h_2 <\ldots<h_k \le h\endSb 1.
$$
Explain why this morally implies (1.10); better still prove (1.10) rigorously (or 
read Gallagher's argument [9]).  
\endproclaim

\proclaim{Exercise 3} We have sketched how the Hardy-Littlewood conjecture implies that 
for a given positive real number $\lambda$, and a fixed non-negative integer $k$, 
 $$
\frac{1}{N} \# \{ n\le N: \ \ \pi(n+\lam \log N)-\pi(n) = k  \} 
 \sim \frac{\lam^k}{k!} e^{-\lam}. 
 $$ 
Deduce that 
$$
\frac{1}{N} \# \Big\{ 2 \le n\le N: \ \ \frac{p_{n+1}-p_n}{\log n} \in [\alpha,\beta]\Big\} 
\sim \int_{\alpha}^{\beta} e^{-t} dt. 
$$ 
\endproclaim 

 \noindent{\bf Proof of (1.10) when $k=2$.}    From the definition (1.5) 
note that ${\frak S}(\{h_1,h_2\}) = {\frak S}(\{0,h_2-h_1\})$ and so, letting $\ell = h_2 -h_1$ 
we see that the LHS of (1.10) is (in the case $k=2$) 
$$
\sum\Sb \ell \le h \endSb {\frak S}(\{0, \ell\}) \Big(\sum\Sb 1\le h_1 < h_2 \le h \\ h_2-h_1=\ell \endSb 
1 \Big) = 
\sum\Sb \ell \le h\endSb {\frak S}(\{0,\ell\}) (h-\ell). 
$$
To evaluate the above asymptotically, it is useful to 
study the generating Dirichlet series
$$
F(s):= \sum_{\ell =1}^{\infty} \frac{{\frak S}(\{0,\ell\})}{\ell^s}. 
$$
The definition (1.5) gives that ${\frak S}(\{0,\ell\}) = \prod_{p|\ell} (1-1/p)^{-1} \prod_{p\nmid \ell} 
(1-2/p)(1-1/p)^{-2}$.  From this, we may see that $F(s)$ converges absolutely in 
the half-plane Re$(s)>1$, and moreover in that region has the Euler product
$$
\align
F(s)&= \prod_{p} \Big( \Big(1-\frac{2}{p}\Big) \Big(1-\frac 1p\Big)^{-2} 
+ \frac{1}{p^s} \Big(1-\frac{1}{p}\Big)^{-1} + \frac{1}{p^{2s}} \Big(1-\frac 1p \Big)^{-1} 
+ \frac{1}{p^{3s}} \Big(1-\frac 1p \Big)^{-1} + \ldots \Big).\\
\endalign
$$
Mutliplying and dividing by $\zeta(s)=\prod_{p} (1-1/p^s)^{-1}$ we see (with a little 
calculation) that in Re$(s) >1$, 
$$
F(s) = \zeta(s) \prod_{p} \Big( 1- \frac{1}{(p-1)^2} + \frac{1}{p^{s-1}(p-1)^2}\Big) = 
\zeta(s) G(s), 
$$
say.  The Euler product for $G(s)$ converges absolutely in Re$(s)>0$ and 
so in that region we have obtained a meromorphic continuation of $F(s)$
with a simple pole at $s=1$ coming from the simple pole of $\zeta(s)$ there.  
We now make use of the formula that for any $c>0$ 
$$
\frac{1}{2\pi i} \int_{c-i\infty}^{c+i\infty} \frac{y^s}{s(s+1)} ds =
\cases 
(1-1/y) &\text{if } y>1\\ 
0 &\text{if } 0< y\le 1.\\
\endcases
$$
This is easily proved by moving the line of integration to the left if $y>1$ 
and to the right if $y\ge 1$; the term $1-1/y$ when $y>1$ arises 
from the residues of the poles at $s=0$ and $s=-1$.  Therefore, if $c>1$,  
we see that 
$$
\align
h \sum\Sb \ell \le h \endSb 
{\frak S}(\{0,\ell\}) \Big(1-\frac{\ell}{h} \Big) 
&= h \sum\Sb \ell =1 \endSb^{\infty} 
{\frak S}(\{0,\ell\}) \frac{1}{2\pi i} \int_{c-i\infty}^{c+i\infty} 
\Big(\frac h\ell\Big)^{s} \frac{ds}{s(s+1)} 
\\
&= \frac{h}{2\pi i} \int_{c-i\infty}^{c+i\infty} 
\zeta(s) G(s) \frac{h^s}{s(s+1)} ds,
\\
\endalign
$$
where the interchange of summation and integration is justified by 
the absolute convergence of $F(s)$ in the region Re$(s)>1$.   To evaluate 
the contour integral, we shift the line of integration to Re$(s)=\epsilon >0$.  In 
the region traversed we encounter only a simple pole at $s=1$ (because of 
$\zeta(s)$) and so our integral is 
$$
h  {\mathop {\text{Res}}_{s=1}} \Big( \frac{\zeta(s)G(s)h^s}{s(s+1)}\Big) + 
\frac{h}{2\pi i} \int_{\epsilon-i\infty}^{\epsilon+i\infty} \zeta(s) G(s) \frac{ds}{s(s+1)}.
$$
Since $G(1)$ is easily seen to be $1$, the residue above equals $G(1)h^2/2 = h^2/2$.
By bounding $\zeta(s)$ and $G(s)$ on the line Re$(s)=\epsilon$ we 
may estimate the remaining integral on that line;  we omit the standard, but technical, details
and merely note that this term is $O(h^{1+\epsilon})$.  We conclude that 
$$
\sum\Sb \ell \le h\endSb {\frak S}(\{0, \ell\}) (h-\ell) 
= \frac{h^2}{2} + O(h^{1+\epsilon}) 
= \sum\Sb h_1 < h_2 \le h \endSb 1 + O(h^{1+\epsilon}).
$$
This proves (1.10) in the case $k=2$.   

 \proclaim{Exercise 4}   Analyze $G(s)$ further by writing it as $\zeta(s+1) H(s)$, 
 where $H(s)$ is now analytic in Re$(s) >- \frac 12$.   Evaluating residues, as above, 
 prove that 
 $$
 \sum_{\ell \le h} {\frak S}(\{0, \ell\}) (h-\ell) = 
 \frac {h^2}{2} - \frac { h\log h}{2} + \frac{Bh} 2+O(h^{\frac 12+\epsilon}),
 $$
 with $B= 1-\gamma-\log 2\pi$; here $\gamma$ is Euler's constant.
\endproclaim 

 \noindent{\bf Concluding remarks.}   Two important consequences of our predictions 
 for the spacings between primes are that 
 $$
 \limsup_{n\to \infty} \frac{p_{n+1}-p_n}{\log n} = \infty, \qquad 
 \text{and} \qquad 
 \liminf_{n\to \infty} \frac{p_{n+1}-p_n}{\log n} = 0.
 $$
 Happily both these results have now been proved.  The first involves 
constructing long strings of composite numbers, and was first proved by 
Westzynthius with important refinements due to Erd{\H o}s and Rankin.  
The second is a recent breakthrough of Goldston, Pintz and \Yildirim, see [12].  
The reader may consult the survey by Heath-Brown [18] for the $\limsup$ result 
and much else besides, and my survey [31] for the $\liminf$ result.  

 \head Lecture 2:  The distribution of primes in longer intervals \endhead 

 \noindent {\bf Cram{\' e}r's prediction.}  In the first lecture we considered the distribution of primes in intervals 
 of length a constant times the average spacing.  We now discuss what happens in longer 
 intervals.  Precisely, we consider $\pi(n+h)-\pi(n)$ for $n\le N$ and where $h/\log N$ is 
 large, but $h/N$ is small.  

 \proclaim{Exercise 5}  Using Stirling's formula, show that as $\lambda$ gets large, 
 a Poisson distribution with parameter $\lambda$ begins to look like a normal 
 distribution with mean $\lambda$ and variance $\lambda$. 
 \endproclaim

 Thus Cram{\' e}r's model would suggest that, if $h/\log N$ is large but $h/N$ is small, 
 then for $n\le N$, $\pi(n+h)-\pi(n)$ has an approximately normal distribution 
 with mean $\sim h/\log N$ and variance $\sim  h/\log N$.  Another way to 
 arrive at this prediction is to calculate the moments (note that for most $n\le N$, 
 $\sum_{\ell=1}^{h} 1/\log (n+\ell) \sim h/\log N$)
 $$
 \frac{1}{N} {\Bbb E} \Big( \sum_{2\le n\le N} \Big( \sum_{\ell= 1}^{h} X(n+\ell) 
 - \sum_{\ell=1}^{h} \frac{1}{\log (n+\ell)} \Big)^k\Big),  
 \tag{2.1a}
 $$
 which we claim is 
 $$
 = \frac{k!}{2^{k/2} (k/2)!} \Big(\frac{h}{\log N}\Big)^{k/2}\Big( 1+ O_k \Big(\frac{\log N}{h} \Big) 
 \Big) \tag{2.1b} 
 $$
 if $k$ is even, and if $k$ is odd it is 
 $$
 \ll \Big(\frac{h}{\log N} \Big)^{(k-1)/2}. \tag{2.1c} 
 $$

 \proclaim{Exercise 6}  Justify (2.1a)-(2.1c) by arguing as follows.  For $n\ge 3$, 
 set $X_0(n) =1-1/\log n$ with probability $1/\log n$ and $X_0(n)= -1/\log n$ with probability 
 $1-1/\log n$: that is, $X_0(n)=X(n)-1/\log n$.   Note that ${\Bbb E}(X_0(n))= 0$.   
 Expand 
 $$
 \frac{1}{N} {\Bbb E} \Big( \sum_{2\le n\le N} \Big(\sum_{1\le \ell \le h} X_0(n+\ell) \Big)^k \Big) 
 = \frac{1}{N} \sum_{\ell_1, \ldots, \ell_k \le h } \sum_{2\le n\le N} 
 {\Bbb E}(X_0(n+\ell_1) \cdots X_0(n+\ell_k) ). 
 $$
 The expectation above is zero if any of the $\ell_i$'s occurs only once among $\ell_1$, $\ldots$, 
 $\ell_k$.  When $k$ is even there is a leading contribution from terms where the $\ell_1$, 
 $\ldots$, $\ell_k$ contain $k/2$ distinct numbers each occurring twice.  
 \endproclaim
 \def\Lam{\Lambda}

\noindent{\bf Calculating the variance via Hardy-Littlewood.}   However, we do 
not believe that this prediction, given by the Cram{\' e}r model, is accurate.  At this juncture, it 
 is more convenient to deal with $\psi(n+h)-\psi(n)$, where $\psi(x)=\sum_{n\le x} \Lambda(n)$
 with $\Lambda(n)$ denoting the von Mangoldt function.  Note that the prime number 
 theorem is equivalent to $\psi(x) \sim x$, and that the Hardy-Littlewood conjecture (1.4) may 
 be recast as (${\Cal H}=\{h_1,\ldots,h_k\}$ is a set of $k$ distinct numbers)
 $$
 \sum\Sb n\le x\endSb \Lam(n+h_1)\cdots \Lam(n+h_k) \sim {\frak S}({\Cal H}) x. \tag{2.2}
 $$
 The Cram{\' e}r model predicts that $\psi(n+h)-\psi(n)$ is approximately normal 
 with mean $\sim h$ and variance $\sim h\log N$.   
 
 To see the flaw in this prediction,
 let us now calculate the variance using the Hardy-Littlewood conjectures.  Note that 
 $$
 \frac{1}{N} \sum_{n\le N} (\psi(n+h)-\psi(n)-h)^2 
 = \frac{1}{N} \sum_{n\le N} \Big(\sum_{\ell \le h} \Lambda(n+\ell)\Big)^2 - 2\frac{h}{N} 
 \sum_{n\le N} \sum_{\ell \le h} \Lam(n+\ell) + h^2.
 $$
 The middle term in the RHS above is $-2h^2 (\psi(N)+O(h\log N))/N\sim -2h^2$.  As for the first 
 term in the RHS we may square it out, and invoke the Hardy-Littlewood conjecture (2.2).
 If we forget all the error terms, then the above is 
 $$
 \frac{1}{N} \sum_{n\le N} \sum_{\ell \le h} \Lam(n+\ell)^2 + 2\sum_{\ell \le h} {\frak S}(\{0, \ell\}) 
 (h-\ell) - h^2.
 $$
 The prime number theorem and partial summation gives that the first term above 
 is $\sim h(\log N -1)$, while from Exercise 4 we see that the second term above is 
 $\sim h^2-h\log h +Bh$.  So, ignoring all error terms, we conclude that the 
 variance satisfies
 $$
 \frac{1}{N} \sum_{n\le N} (\psi(n+h)-\psi(n)-h)^2 \sim h \Big(\log \frac Nh + B-1\Big), \tag{2.3}
 $$ 
which is different from the $\sim h\log N$ predicted by Cram{\' e}r's model. 

 \proclaim{Exercise 7}  Assume that the Hardy-Littlewood conjecture (2.2) 
 holds in the quantitative form 
 $$
 \sum_{n\le x} \Lam(n+h_1) \cdots \Lam(n+h_k) = {\frak S}({\Cal H}) x +O(x^{\frac 12+\epsilon}),
 $$
 uniformly for $k\le K$, and distinct $h_j$ satisfying $1\le h_j \le x$.  Using this, obtain (2.3) 
 with an error term of $O(h^{\frac 12+\epsilon} + h^2 N^{-\frac 12+\epsilon}+h^3N^{-1})$.  
 Thus, even assuming the quantitative Hardy-Littlewood conjectures, one knows 
 (2.3) only for $h\le N^{\frac 12 -\epsilon}$.
 \endproclaim

 So although the Hardy-Littlewood probabilities and the Cram{\' e}r probabilities 
 are roughly equal on average, significant deviations show up when we consider $h$ 
 to be a small power of $N$.  We believe that (2.3) is the right asymptotic for the 
 variance and the Cram{\' e}r model predicts the wrong answer.   

 \noindent {\bf The variance and zeros of the zeta function.}  Here is the sketch of a 
 very different calculation which leads to the same answer as (2.3).    Riemann's explicit 
 formula (see [6]) says that 
 $$
 \psi(x) = x -\sum_{\rho} \frac{x^{\rho}}{\rho} + \text{negligible terms}.
 $$
 Here $\rho$ runs over the non-trivial zeros of the Riemann zeta-function.
We assume the Riemann hypothesis and write $\rho=\tfrac 12+i\gamma$.  The sum over 
zeros is only conditionally convergent, but we will argue loosely omitting such considerations.  
It follows that 
$$
\psi(x+h)-\psi(x)-h = -\sum_{\rho} \frac{(x+h)^{\rho}-x^{\rho}}{\rho}  + \text{negligible terms}.
$$
The sum over zeros above is weighted down with a factor $1/\rho$, and so 
we may expect that large zeros make a minor contribution.  It turns out that 
zeros with $|\rho| \ge x/h$ are not so important.  For the small zeros, we 
replace $(x+h)^{\rho}-x^{\rho}$ by the Taylor approximation $\rho h x^{\rho-1}$.  
Therefore we may expect that 
$$
\align
\frac 1X \int_X^{2X} (\psi(x+h)-\psi(x)-h)^2 dx &\approx 
\frac{h^2}{X^2} \int_{X}^{2X} \Big|\sum_{|\gamma|\le X/h} x^{i\gamma} \Big|^2 dx 
\\
&= \frac{h^2}{X} \sum_{|\gamma_1|, |\gamma_2| \le X/h} 
X^{i(\gamma_1-\gamma_2)} \frac{2^{1+i(\gamma_1-\gamma_2)}-1}{1+i(\gamma_1-\gamma_2)}.
\tag{2.4}\\
\endalign
$$
There are $\ll \log T$ zeros of the zeta-function with ordinates lying between $T$ and $T+1$.  
Using this observation in (2.4), and estimating the magnitude of the sums over zeros there, 
 we ``deduce" that, assuming RH, \footnote{In fact we would only deduce $\ll h(1+\log X/h)^3$ 
 but the extra ``$\log$" may be removed by smoothing.}
$$
\frac 1X\int_X^{2X} (\psi(x+h)-\psi(x)-h)^2 dx \ll h (1+\log X/h)^2.  \tag{2.5}
 $$
 A result like this was established by  Selberg [30].    If we want an asymptotic in (2.4), 
 then we need some understanding of the spacings $\gamma_1-\gamma_2$ between zeros of the Riemann  zeta-function.  Such an understanding is provided by the pair correlation conjecture 
 of Montgomery [24], which predicts that these ordinates are distributed like eigenvalues of large 
 random matrices.  Using such information Mueller obtained an 
 asymptotic formula much like (2.3), and Goldston and Montgomery [11] showed conversely 
 that a formula like (2.3) also conveys information about the zeros of $\zeta(s)$.  For more discussion 
 on this set of ideas consult Goldston's recent survey [10].

\noindent {\bf Higher moments.}  Recently, Montgomery and I (see [25] and [26]) used a quantitative form of the Hardy-Littlewood 
 conjecture (see Exercise 7) to study higher moments of $\psi(n+h)-\psi(n)-h$.  We now describe 
 these results briefly.  They support the conjecture that if $(\log N)^{1+\delta} 
 \le h\le N^{1-\delta}$ then for $n\le N$ the distribution of $\psi(n+h) 
 -\psi(n)$ is approximately normal with mean $h$ and variance $h\log (N/h)$.

 We assume that $(\log N)^{1+\delta} \le h\le N^{1-\delta}$ and wish to evaluate 
 $$
 \frac 1N \sum_{n\le N} (\psi(n+h)-\psi(n)-h)^r. \tag{2.6} 
 $$
For even $r$ we expect that this is $\sim \frac{r!}{2^{r/2}(r/2)!} (h\log N/h)^{r/2}$, 
while for odd $r$ we expect it to be $o((h\log N/h)^{r/2})$.   If we simply expanded 
$(\psi(n+h)-\psi(n)-h)^{r}$ in powers of $(\psi(n+h)-\psi(n))$ and $h$ (as we did in 
the case $r=2$) then we would get many terms all of size $h^r$, and a careful 
cancellation of these and lower order terms is needed before we get to the actual delicate 
main term of size essentially $h^{r/2}$.  To circumvent this, we define $\Lam_0(n) 
=\Lam(n)-1$, in analogy with Exercise 6.  This eliminates the unnecessary higher order terms at 
the outset, and simplifies calculations considerably.  For other situations where 
this trick helps, see my paper with Granville [15] in this volume.  
Using this notation, and expanding (2.6) we want to understand
 $$
\sum\Sb \ell_1, \ldots, \ell_r \le h \endSb  \frac{1}{N} \sum_{n\le N} \Lam_0(n+\ell_1)\cdots \Lam_0(n+\ell_r). \tag{2.7}
$$

 \proclaim{Exercise 8}  Define the modified singular series ${\frak S}_0({\Cal H})$ 
 by 
 $$
 {\frak S}_0({\Cal H}) = \sum_{{\Cal J} \subset {\Cal H}} (-1)^{|{\Cal H}| - |{\Cal J}|} {\frak S}({\Cal J }), 
 \qquad \text{so that } \ \ {\frak S}({\Cal H}) =\sum_{{\Cal J}\subset {\Cal H}} {\frak S}_0({\Cal J}).
 $$
 Here we understand that ${\frak S}(\emptyset) = {\frak S}_0(\emptyset) =1$.  Show that 
 the quantitative Hardy-Littlewood conjecture of Exercise 7 is the same as
 $$
 \sum_{n\le x} \Lam_0(n+h_1)\cdots \Lam_0(n+h_k) = {\frak S}_0({\Cal H}) x +O(x^{\frac 12+\epsilon}),
 $$
keeping the hypotheses there.
 \endproclaim

 For simplicity, consider first the terms in (2.7) when the $\ell_i$ are distinct.   If we use the 
 asymptotic of Exercise 8 we are led to the problem of evaluating 
 $$
 \sum\Sb h_1, \ldots , h_k \le h \\ h_i \text{ distinct} \endSb {\frak S_0}({\Cal H}),  
 $$
 which is a problem analogous to, but more delicate than, Gallagher's calculation (1.10).  
 The main result in [26] is the asymptotic 
 $$
  \sum\Sb h_1, \ldots , h_k \le h \\ h_i \text{ distinct} \endSb {\frak S_0}({\Cal H}) \ = \qquad
   \cases 
  \{ 1+o(1)\} \  \frac{k!}{2^{k/2}(k/2)!} (-h\log h + B+1)^{k/2} &\text{if $k$ is even}\\ 
  o((h\log h)^{k/2}) &\text{if $k$ is odd}. \\
  \endcases
  \tag{2.8}
  $$

 \proclaim{Exercise 9}  Show the following refinement of Gallagher's (1.10):  
 $$
\sum\Sb h_1, \ldots, h_k \le h\\ h_i \text{ distinct} \endSb 
{\frak S}({\Cal H}) 
= h^k - \binom{k}{2} h^{k-1} \log h + \binom{k}{2} B h^{k-1} + O(h^{k-3/2+\epsilon}). 
$$
\endproclaim

 Returning to (2.7), we must analyze the terms when the $\ell_i$ are not necessarily 
 distinct.   Suppose that $h_1$, $\ldots$, $h_k$ are the distinct elements among 
 $\ell_1$, $\ldots$, $\ell_r$ and that each $h_i$ appears $m_i (\ge 1)$ times 
 among the $\ell_i$.  After a little combinatorics, we may write (2.7) as 
 $$
 \sum_{k=1}^{r} \sum\Sb m_1, \ldots, m_k \ge 1 \\ \sum m_i =r \endSb 
 \binom{r}{m_1,\ldots, m_k} \frac{1}{k!} \sum\Sb h_1, \ldots, h_k \le h\\ 
 h_j \text{ distinct} \endSb \frac{1}{N} \sum_{n=1}^{N} \prod_{i=1}^{k} \Lam_0(n+h_i)^{m_i}. 
 \tag{2.9}
 $$
 We must distinguish the indices where $m_i=1$ and the remaining indices where $m_i>1$.  
 Let ${\Cal I}$ denote the subset of $\{1,\ldots,k\}$ such that $m_i=1$ for $i\in {\Cal I}$.  
 For $i\notin {\Cal I}$ (so $m_i \ge 2$) 
 we think of $\Lam_0(n+h_i)^{m_i}$ as being essentially $(\log N)^{m_i-1} \Lam(n+h_i)$: the point is 
 that both quantities have about the same expected value $(\log N)^{m_i-1}$, 
 unlike the case when $m_i=1$ where the expected value of $\Lam(n+h_i)$ and $\Lam_0(n+h_i)$ are $1$ and $0$ respectively.  Therefore the inner sum over $n$ in (2.9) is essentially 
 $$
 \align
&\frac{(\log N)^{r-k}}{N} \sum_{n=1}^{N} \prod_{i\in {\Cal I}} \Lam_0(n+h_i) 
\prod\Sb 1\le i\le k\\ i\notin {\Cal I}\endSb (\Lam_0(n+h_i)+1) 
\\
=&\frac{(\log N)^{r-k}}{N} \sum\Sb {\Cal I} \subset {\Cal J} \subset \{1,\ldots,k\} \endSb 
\sum_{n=1}^{N} \prod_{j\in {\Cal J}} \Lam_0(n+h_j).
\\
\endalign
$$
Now we invoke the Hardy-Littlewood conjecture of Exercise 8, and use (2.8). 
 
\proclaim{Exercise 10}  Complete the details in evaluating (2.7).  Show that when $r$ is odd 
or any of the $m_i$'s is $\ge 3$ we get a contribution of $o(h\log N)^{r/2}$.   In the case $r$ is 
even, the main term $\frac{r!}{2^{r/2}(r/2)! }(h\log N/h)^{r/2}$ arises from contributions to (2.9) 
where the $m_i$ are all $1$ or $2$.  
\endproclaim

 The proof of (2.8) is quite complicated, and we do not go into it here.  Let us 
 however point out one important ingredient.  While motivating 
 the Hardy-Littlewood conjecture in Lecture 1, we considered the toy problem of 
 reduced residues $\pmod q$.   If  $1=a_1 < a_2 <\ldots < a_{\phi(q)} < q$ are the reduced residues below $q$, then  we may ask for the distribution of $(a_{i+1}-a_i)(\phi(q)/q)$; we have 
 multiplied by $\phi(q)/q$ so that this is $1$ `on average.'  If, for example, $q$ is the product of 
 the first $\ell$ primes, then as in Lecture 1 we may think of these $a_i$ as being like primes, 
 and expect that, for $0<\alpha < \beta$, 
 $$
 \# \Big\{ 1\le i \le \phi(q) -1 : \ \ (a_{i+1}-a_i) \frac{\phi(q)}{q} \in [\alpha,\beta] \Big\}  
 \sim \phi(q) \int_{\alpha}^{\beta} e^{-x} dx.
 $$
 A beautiful result of Hooley [20] shows that this holds provided $\phi(q)/q$ is small.  Obviously, some 
 restriction on $\phi(q)/q$ is necessary; for example, if $q$ is prime then clearly $a_{i+1}-a_i=1$ 
 for $1\le i\le \phi(q)-1$.   Moreover, Montgomery and Vaughan [27] have even estimated the 
 moments: 
 $$
 M_k(q;h)=\sum_{n=1}^{q} \Big(\sum\Sb \ell\le h\\  (n+\ell,q)=1\endSb 1-h\frac{\phi(q)}{q}\Big)^k.
 $$
 The proof of (2.8) builds on the techniques developed there.  

 In our discussion above we have ignored error terms altogether.  If one argues carefully 
 using the quantitative Hardy-Littlewood conjectures of Exercises 7 and 8, we can 
 evaluate the $r$-th moment (2.6) provided that $h\le N^{1/r-\epsilon}$.  We 
 expect that the same asymptotics hold even when $h$ is larger with $h\le N^{1-\epsilon}$.  
 Thus these arguments suggest that for $(\log N)^{1+\delta} \le h\le N^{1-\delta}$, the 
 distribution of $\psi(n+h)-\psi(n)$ (for $n\le N$ is approximately normal with mean $h$ and
 variance $h\log N/h$.  For numerical support for this conjecture, see [3] and [25].  For other work related
 to this circle of ideas, see [3] and [4].

 \noindent{\bf Connections with zeros of $\zeta(s)$?}  We mentioned earlier the work of 
 Goldston and Montgomery relating the variance of primes in short intervals to the pair correlation 
 of zeros of $\zeta(s)$.    Our calculations on the higher moments of primes in short intervals 
suggest that if $X \ge T^{1+\epsilon}$ then\footnote{We take this opportunity to point out 
that the important constraint $X\ge T^{1+\epsilon}$ has been erroneously omitted in 
a similar discussion (on page 594) in [26].}
$$
\int_X^{2X} \Big(\sum_{|\gamma| \le T} x^{i\gamma}\Big)^k dx =\int_X^{2X} 
\Big(\sum_{0\le \gamma \le T} 2 \cos(\gamma \log x)\Big)^{k} dx
$$
is $\sim \frac{k!}{2^{k/2}(k/2)!} X (2N(T))^{k/2}$ if $k$ is even, and is $o(XN(T)^{k/2})$ if $k$ is 
odd.  Here $N(T)\sim \frac{T}{2\pi} \log T$ denotes the number of zeros of 
$\zeta(s)$ with $0\le \gamma \le T$.  Viewed this way, Montgomery's pair correlation 
conjecture may be thought of as saying that for  $x\ge T^{1+\epsilon}$ the sum $\sum_{0\le 
\gamma\le T} \cos(\gamma \log x)$ behaves like a sum of uncorrelated random variables. 
The higher moments suggest that it behaves in fact like a sum of independent random 
variables. \footnote{This is analogous to a result of E. Rains [28] in random matrix theory.} 
 These statements are quite vague, and it would be nice to flesh out the 
precise connection between these higher moments and zeros of $\zeta(s)$.  
For other connections between zeros of $\zeta(s)$ and Hardy-Littlewood type 
conjectures see [2].

 \noindent{\bf Chebyshev's bias.}  We have considered above the distribution of primes 
 in short intervals.  What happens to the distribution in long intervals $[1,x]$?   That is 
 what can be said about the distribution of $\psi(x)-x$.  Assuming RH, we 
 get from Riemann's explicit formula that this is essentially 
 $-2 x^{\frac 12}\text{Re } \sum_{0<\gamma}  x^{i\gamma}/(1/2+i\gamma)$.  It is 
 expected\footnote{There is perhaps no good reason for this belief, except that 
 the contrary situation is harder to imagine!}
 that the zeros of $\zeta(s)$ are all simple, and have no non-trivial ${\Bbb Q}$-linear relations 
 among them.  In that case the sum over zeros above may be modeled by $\text{Re }
 \sum_{0<\gamma} X(\gamma)/(1/2+i\gamma)$ where the $X(\gamma)$ are independent 
 random variables, taking uniformly distributed values on the unit circle.  Precisely, as 
 $t$ varies from $1$ to $T$, the distribution of $(\psi(e^t)-e^t)/(2e^{t/2})$ is 
 like the distribution of our random sum above.\footnote {The change of variable $x=e^{t}$ 
 means that $x^{i\gamma}=e^{it\gamma}$ now takes values uniformly on the unit 
 circle as $t$ varies.}  This is a certain non-universal distribution, which has been 
 investigated in, for example, [23] and [29].  To gain a flavor of this distribution the reader 
 may contemplate $\sum_{n=1}^{\infty} X_n/n$ where the $X_n$ are independent 
 random variables taking the values $\pm 1$ with equal probability.   

 The distribution above is symmetric about the origin, and so $\psi(x)$ is as likely to be 
 larger than $x$ as it is to be smaller than $x$.  However, $\psi(x)=\theta(x)+\theta(x^{\frac 12}) +
 \theta(x^{\frac 13})+\ldots$ where $\theta(x)=\sum_{p\le x} \log p$.  Thus it 
 is much more likely for $\theta(x)$ to be smaller than $x$ than for it to be larger than $x$.  
 By partial summation one gets that $\pi(x)< \text{li}(x)$ much more often than $\pi(x)>\text{li}(x)$.  
 In fact, in a certain sense the probability that $\pi(x)$ `beats' $\text{li}(x)$ is only  
 $0.00000026\ldots$!  We stop here, referring the reader to Rubinstein and Sarnak [29], and 
 the delightful recent survey [14] for more information.  

 To summarize, we found three distinct behaviors for the distribution of primes in intervals.  
 At the ``microscopic" scale ($h \asymp \log N$) there is Poisson behavior, 
 at the ``mesoscopic" scale ($h/\log N\to \infty$, $h=o(N)$) there is Gaussian (normal) behavior, 
 and at the ``macroscopic" scale ($h \gg N$) there is a specific non-universal distribution law.  
 Such division into three regimes occurs in many other problems as well;  for example, 
 in the distribution of lattice points in the plane.  As a starting point, we refer the reader 
 to the recent paper of Hughes and Rudnick [21] and to the references therein.

 \head Lecture 3:  Maier's method and an ``uncertainty principle" \endhead

 \noindent If the Riemann Hypothesis is true, then from Selberg's result (2.5) we easily deduce that (for $h\le N$) 
 the number of $n\le N$ with $|\psi(n+h)-\psi(n)-h| \ge \sqrt{h} (\log N)^{1+\delta}$ 
 is $\ll N/(\log N)^{2\delta}$.   It follows that if $N\ge h\ge (\log N)^{2+\delta}$ then 
 ``almost all" intervals $(n,n+h]$ with $n\le N$ contain about the correct number of primes,
$\sim h/\log N$.  
If (2.3) holds then we can even conclude that if $h\ge (\log N)^{1+\delta}$ then 
``almost all" intervals $(n,n+h]$ with $n\le N$ contain approximately the correct number of 
primes.   In Cram{\' e}r's model, one can show that almost surely $\sum_{x\le n\le x+h} X(n)
\sim h/\log x$ if $(\log x)^{2+\delta}\le h\le x$.  Thus it seems quite plausible that 
if $x$ is large and $x\ge h\ge (\log x)^{2+\delta}$ then $\psi(x+h)-\psi(x)\sim h$.

The classical prime number theorem with error term $x\exp(-C\sqrt{\log x})$ tells us that such 
a result holds if $h\ge x\exp(-C\sqrt{\log x})$.   An important advance was made by  Hoheisel 
who showed that the asymptotic $\psi(x+h)-\psi(x)\sim h$ holds if $x\ge h\ge x^{\theta}$ 
for some number $\theta<1$.  He was able to take $\theta=1-1/33000$, but this 
has been improved subsequently, with the best result known, due to Huxley, being $\theta=
\frac 7{12}+\epsilon$.   If the Riemann hypothesis is true then $\theta$ may be taken as $\frac 12+\epsilon$.  The arguments pioneered by Hoheisel depend on the fact that while we don't know 
RH, we do know that  most zeros of $\zeta(s)$ lie close to the $\frac 12$ line.  For a nice account of 
these results see Heath-Brown [18].   If the asymptotics given for (2.6) are true then 
we may take $\theta$ to be any positive number.

Thus it seems the conjecture that $\psi(x+h)-\psi(x) \sim h$ for $x\ge h \ge (\log x)^{2+\delta}$,
if true, lies 
quite deep.  This conjecture was widely believed until the mid 1980's when Maier [22] shattered 
this belief by showing that for any $A>1$ there are arbitrarily large $x$ such 
that the interval $(x,x+(\log x)^A]$ contains significantly more primes than usual 
(that is, $\ge (1+\delta_A) (\log x)^{A-1}$ primes for some $\delta_A>0$) and 
also intervals $(x,x+(\log x)^A]$ containing significantly fewer primes than usual.  In this lecture 
we will sketch Maier's ingenious method, and describe some extensions of his idea.  The 
reader may also consult [13] for another exposition of related ideas.

\noindent {\bf Maier's ``matrix" method.}  Let $x$ be large, and $h$ be on the 
scale of a power of $\log x$.  Let $P$ be an integer which we will 
eventually take to be the product of many small (below $\log x$) primes.  
Consider the $[x/P]$-by-$h$ ``matrix" whose $(i,j)$-th entry is the number 
$([x/P]+i)P+j$.  Thus the entries of this matrix are numbers lying between $x$ and $2x+h$.  
Note that each row of the matrix consists of an interval $([x/P]+i)P$ to $([x/P]+i)P +h$.  
Each column of the matrix consists of an arithmetic progression with common difference $P$: 
namely, $x\le n\le 2x+h$ with $n\equiv j\pmod P$.   The idea is to count the number of primes in 
this matrix in two ways: by counting the primes row by row, and by  counting the 
primes column by column, and then comparing the two answers.  If we assume that the asymptotic formula for primes in short intervals 
holds then we get an answer for the row by row calculation.  The prime number theorem 
for arithmetic progressions allows us to do the column by column calculation.  Of course the two 
answers should match.  However when $h$ is very small, like a power of $\log x$, there are 
choices of $P$ for which the answers don't match!  This leads to Maier's result.

Consider the row by row calculation.  The number of primes is 
$$
\sum\Sb x/P \le n \le 2x/P \endSb (\pi(nP+h)-\pi(nP)), \tag{3.1a}
$$ 
and if we assume that intervals of length $h$ contain the right number of primes, this is 
$$
\sim \frac{x}{P} \frac{h}{\log x}.  \tag{3.1b}
$$

Consider next the column by column calculation.  If the progression $n\equiv j\pmod P$ 
is to contain primes, we must have $(j,P)=1$.  In that case the prime number 
theorem in arithmetic progressions would say that such a progression contains 
a proportion $1/\phi(P)$ of all primes.  Of course, in order to use the prime number theorem in 
arithmetic progressions rigorously we must pay attention to the size of the modulus $P$ compared 
with $x$.  Assuming that this is not an issue, we find that the column by column contribution is 
$$
\sum\Sb j\le h\\ (j,P) = 1\endSb (\pi(2x+h;P,j)-\pi(x;P,j)) 
\sim \frac{x}{\phi(P)\log x} \sum\Sb j\le h\\ (j,P)=1\endSb 1.  \tag{3.2}
$$

If we compare (3.2) and (3.1b) we find the relation 
$$
\sum\Sb j\le h \\ (j,P)=1 \endSb 1 \sim h \frac{\phi(P)}{P}  \tag{3.3}
$$
should hold.
At first glance, (3.3) is eminently reasonable: the probability that $j$ is coprime to $P$ is 
$\phi(P)/P$.  It is even easy to make this precise: write the condition $(j,P)=1$ as $\sum_{\ell | (j,P)} \mu(\ell)$ and we easily get 
$$
\sum\Sb j\le h \\(j,P) =1 \endSb 1 = h\frac{\phi(P)}{P} + O(d(P)), \tag{3.4} 
$$
where $d(P)$ is the number of divisors of $P$.  Thus, if $h$ is just a bit larger than $d(P)$ 
(which is always quite small, that is $\ll P^{\epsilon}$) then (3.3) will hold.   So where is the contradiction?  
The point is that in Maier's application $h$ is very small compared with $P$, and so (3.4) 
is useless.  

For the purpose of illustration suppose that $P$ is the product of all
primes between $(\log x)^{\frac{9}{10}}$ and $(\log x)/100$.  Then, by the prime number 
theorem, $P$ is about size $x^{\frac{1}{100}+o(1)}$.  For such moduli $P$ we 
don't know the prime number theorem in arithmetic progressions used in (3.2), 
but such a result does hold if the Riemann hypothesis for Dirichlet $L$-functions is 
true; let us postpone a discussion of this point.  Suppose now that $h$ is a number of size $(\log x)^{\theta}$ with $2 < \theta < 2.7$.  
By inclusion-exclusion, the LHS of (3.3) is 
$$
\align
&\sim h -\sum_{(\log x)^{0.9} \le p\le (\log x)/100} \frac{h}{p} 
+ \sum\Sb (\log x)^{0.9} \le p < q \le (\log x)/100 \endSb \frac{h}{pq} 
\\
&\sim h \Big( 1- \log \frac{10}{9} + \frac{1}{2} \Big(\log \frac{10}{9}\Big)^2 \Big),
\\
\endalign
$$
where we have used the prime number theorem to evaluate $\sum 1/p$ for $p$ between 
$(\log x)^{\frac 9{10}}$ and $(\log x)/100$.  On the other hand, by Mertens' theorem,
the RHS of (3.3) is 
$$
\sim h \prod_{(\log x)^{\frac{9}{10}} \le p\le (\log x)/100} \Big(1-\frac{1}{p}\Big)^{-1} 
\sim \frac{9}{10} h.
$$
 The formula for the LHS has the first three terms in the usual expansion of 
 $\frac 9{10} = e^{-\log (10/9)}$, so the two answers are certainly close, 
 but obviously they are not equal!   Indeed the LHS is a little bit larger.  

 \proclaim{Exercise 11}  Conclude from the above that for any $2 < \theta <2.7$ there 
 exist arbitrarily large $x$ such that the interval $[x,x+(\log x)^{\theta}]$ contains 
 significantly more primes than expected.  Taking such an interval and cutting 
 it up into smaller intervals,   deduce that the same conclusion holds for all $1< \theta <2.7$.  
 Using the same $P$ as above, and taking four terms in the inclusion-exclusion 
 formula, show that if $\theta<3.6$ there exist intervals $[x,x+(\log x)^{\theta}]$ with 
 significantly fewer primes than expected.   In this manner one can proceed for $\theta< 8.1$, just using 
 inclusion-exclusion and easy calculations.   Now replace $(\log x)^{0.9}$ in 
 the definition of $P$ with $(\log x)^{1-\delta}$ and prove  Maier's theorem.  
 \endproclaim

 \noindent{\bf More on the contradiction.}   Now let us describe a different way of 
 seeing a contradiction to (3.3).  This method is very flexible, and works for many choices of $P$, 
 and also generalizes readily.  Let $y$ be a large parameter; in the application we 
 may think of $y$ as being some power of $\log x$.  From each dyadic block 
 $[2^{-j}y, 2^{-j+1}y]$ with $j\le [\log y/(2\log 2)]$ select about half the primes.  
 Take $P$ to be the product of these selected primes.  Thus $P$ is composed of about half the 
 primes in $[\sqrt{y},y]$, and there are plenty of choices for $P$.  
 Let $u\ge 1$ be a real number, set $h=y^{u}$ 
and consider whether (3.3) can hold.    We will show that for arbitrarily large $u$ the 
LHS is appreciably larger than the RHS, and for arbitrarily large $u$ it is smaller.

 To see this we consider the Dirichlet series $\zeta_P(s) =\sum_{(n,P)=1} n^{-s}$.   Plainly we 
 have 
 $$
 \zeta_P(s) = \zeta(s) \prod_{p|P} \Big(1-\frac 1{p^s}\Big), \tag{3.5}
 $$
 so that $\zeta_P(s)$ extends to a meromorphic function in all of ${\Bbb C}$ with a 
 simple pole at $s=1$.  The point is that if something like (3.3) holds then $\zeta_P(s)$ 
 must approximately look like $\zeta(s)\phi(P)/P$, and by choosing $s$ appropriately we 
 can obtain a contradiction to (3.5).  More precisely, set 
 $$
 E(u) = \frac{1}{y^u} \Big( \sum\Sb n\le y^u \\ (n,P)=1\endSb 1 -[y^u] \frac{\phi(P)}{P}\Big).  
 $$
 Then,  for Re$(s)>1$, 
 $$
 \zeta_P(s)-\zeta(s)\frac{\phi(P)}{P}  = \int_{1^-}^{\infty} \frac{1}{z^s} d\Big(\sum\Sb n\le z\\ (n,P)=1 
 \endSb 1 - \sum_{n\le z} \frac{\phi(P)}{P}\Big) 
 = \int_1^{\infty} \frac{s}{z^{s}} E\Big(\frac{\log z}{\log y}\Big) dz,
$$
upon integrating by parts.  Changing variables $u=\log z/\log y$ we obtain that 
$$
\zeta_P(s) = \zeta(s) \frac{\phi(P)}{P} + s\log y \int_0^{\infty} E(u) y^{-u(s-1)} du. \tag{3.6}
$$
To start with, (3.6) is valid for Re$(s)>1$, but since $E(u) \ll d(P)y^{-u}$ by (3.4), we 
see that (3.6) makes sense for Re$s>0$.  

 \proclaim{Exercise 12}  Let $(\log y)/2 \ge \xi \ge 1$ be a real number, and 
 take $s=1-\xi/\log y + i\pi/\log y$.  Using (3.5) prove that 
 $$
 |\zeta_P(s)| \gg \frac{\log y}{\xi} \exp\Big(   \frac{e^{\xi}}{2\xi} + O\Big(\frac{e^{\xi}}{\xi^2}\Big)\Big).
 $$
 Then using (3.6) deduce that 
 $$
 \int_0^{\infty} |E(u)| e^{\xi u} du \gg \frac{1}{\xi}  \exp\Big(   \frac{e^{\xi}}{2\xi} + O\Big(\frac{e^{\xi}}{\xi^2}\Big)\Big).
 $$
 Show that $\int_0^{\infty} E(u)e^{\xi u} du \ll 1/\xi$, so that in the LHS above both 
 positive and negative values of $E(u)$ make roughly equal contributions. 
 \endproclaim 

 Morally, Exercise 12 shows that $E(u)$ cannot be too small for large $u$.  To make this 
 precise, one also needs an upper bound for $E(u)$ so as to be able to bound the 
 tail of the integral in Exercise 12.  Developing this argument carefully, one may show 
 that there is a positive constant $A$ such that every interval $[u(1-A/\log u),u(1+A/\log u)]$ 
 contains points $u_\pm$ satisfying 
 $$
 E(u_+) \ge \exp(-u_+ (\log u_+ + \log \log u_+ + O(1))),
 $$
 and 
 $$
 E(u_-) \le - \exp(-u_-(\log u_- + \log \log u_- + O(1))).
 $$
 For more details, see section 3 of [16], especially Corollary 3.3.

 Earlier, we postponed discussion of the prime number theorem in arithmetic progressions.  
 We refer the reader to Davenport [6] for an account of this.  In Chapter 20 there one finds 
 Page's result that $\psi(x;q,a) \sim x/\phi(q)$ for all $q\le \exp(C\sqrt{\log x})$ with 
 the possible exception of multiples of a particular modulus $q_1$ which may depend on $x$.  
 If we choose $y$ a little less than $\sqrt{\log x}$ then our moduli $P$ above are below $\exp(C\sqrt{\log x})$ and certainly we can find $P$ that are not multiples of the exceptional modulus $q_1$.  
 Thus our appeal to the prime number theorem in arithmetic progressions can be made 
 rigorous.  

 The flexibility in choosing $P$ is quite useful.  Exploiting this, Granville and I (see [16]) showed that 
 the asymptotic 
 $$
 \psi(x+h)-\psi(x) = h +O(h^{\frac 12+\epsilon}), \tag{3.7} 
 $$
suggested by Cram{\' e}r's model, sometimes fails to hold if $h\le \exp((\log x)^{\frac 12-\epsilon})$.  
This improves work of Hildebrand and Maier [19] who had obtained this result assuming the 
Generalized Riemann Hypothesis, and a weaker result unconditionally.  It seems safe to conjecture 
that (3.7) holds if $h \ge x^{\delta}$, and perhaps it holds when $h\ge \exp((\log x)^{\frac 12+\delta})$.

 \noindent {\bf An uncertainty principle.}   Maier's method can be adapted to establish 
 limitations to the equidistribution of primes in arithmetic progressions.  For example, 
 Friedlander and Granville [8] proved that for every $A\ge 1$ there exist 
 large $x$ and an arithmetic progression $a\pmod q$ with $(a,q)=1$ and $q\le x/(\log x)^A$ 
 such that 
 $$
 \Big| \pi(x;q,a) -\frac{\pi(x)}{\phi(q)}\Big| \gg_A \frac{\pi(x)}{\phi(q)}.
 $$

 More recently, Balog and Wooley [1] showed that the sequence of integers which 
 may be written as the sum of two squares also exhibits ``Maier type" irregularities 
 in intervals $(x,x+(\log x)^A)$ for any fixed, positive $A$.  Previously Maier's work 
 had seemed inextricably linked to the mysteries of primes, but Balog and Wooley's 
 result suggests that such results should be part of a more general phenomenon.  
 This has been formalized by Granville and me as an ``uncertainty principle" for arithmetic 
 sequences.  What Maier's argument shows is that the primes cannot be simultaneously 
 well distributed in short intervals, and in arithmetic progressions.  Then a suitable version 
 of the prime number theorem in arithmetic progressions is used to remove the second 
 possibility, leaving us with the irregularities of distribution in short intervals.   The first conclusion 
 of irregularities in short intervals or progressions turns out to be a general feature of 
 many interesting arithmetical sequences.  

A rough description of this result is as follows:  Let ${\Cal A}$ denote a sequence $a(n)$ of non-negative 
 real numbers, and let ${\Cal A}(x) =\sum_{n\le x}a(n)$.  If ${\Cal A}$ is well-distributed in 
 short intervals, then we may expect that 
 $$
 {\Cal A}(x+y) - {\Cal A}(x) \approx y \frac{{\Cal A}(x)}{x} .
  \tag{3.8}
 $$
 To understand the distribution of $a(n)$ in arithmetic progressions we begin with 
 $n$ that are multiples of a given number $d$.  We suppose that there is a non-negative 
 multiplicative function $h$ such that 
 $$
 {\Cal A}_d(x) = \sum\Sb n\le x\\ d|n \endSb a(n) \approx \frac{h(d)}{d} {\Cal A}(x). \tag{3.9a}
 $$
 We assume that the asymptotic behavior of 
 $$
 {\Cal A}(x;q,a):=\sum\Sb n\le x \\ n\equiv a\pmod q\endSb a(n) 
 $$
 depends only on the g.c.d. of $a$ and $q$.  Then (3.9a) leads to the 
 prediction that 
 $$
 {\Cal A}(x;q,a) \approx \frac{f_q(a)}{q\gamma_q} {\Cal A}(x), \tag{3.9b}
 $$
 with $\gamma_q =\prod_{p|q} (p-1)/(p-h(p))$ and $f_q(a)$ is a certain non-negative 
 multiplicative function of $a$, defined in terms of $h$, such that $f_q(a)=f_q((a,q))$ 
 so that $f_q(a)$ is periodic $(\bmod\  q)$.   We can be flexible in how we want 
 to  assume (3.9b); for example,  sometimes it is convenient to assume it 
 only for $q$ that are coprime to a certain fixed set of primes.  

 To illustrate the framework consider the following examples.  

 \noindent {\bf Example 1.}  Take $a(n)=1$ for all $n$.  It is natural to take $h(d)=1$ 
 for all $d$, $\gamma_q=1$, and $f_q(a)=1$.  Then (3.9a) and (3.9b) are both 
 good approximations with errors at most $1$. 

\noindent{\bf Example 2}.  Take $a(n)$ to be the indicator function of the primes.  
Then $h(1)=1$ and $h(d)=0$ for $d>1$.  One has $\gamma_q= \phi(q)/q$ and 
$f_q(a)=1$ if $(a,q)=1$ and $0$ otherwise.  The prime number 
theorem in arithmetic progressions gives (3.9b) for small values of $q$.  
The result of Friedlander and Granville places restrictions on the approximation (3.9b) when 
$q$ is large.  Maier's results place restrictions on (3.9a) for small $y$.  

\noindent{\bf Example 3}.  Take $a(n)$ to be the indicator function of the sums 
of two squares.   The multiplicative function $h$ is defined by $h(p^k)=1$ if $p^k\equiv 1\pmod 4$ 
and $h(p^k)=1/p$ if $p^k \equiv 3\pmod 4$.   Here Balog and Wooley's result places 
restrictions on (3.9a).

The main results of [16] give that  if $h(p)$ is not always close to $1$ (as in the regular 
example 1) then there will be moduli $q$ for which (3.9b) cannot hold.  Typically 
these moduli will be large as in the Friedlander-Granville result for primes in 
progressions.  Furthermore, either there exist values $y$ larger than an arbitrary 
power of $\log x$ for which (3.9a) is false, or there exist small moduli $q$ 
(below $\exp((\log x)^{\delta})$) for which (3.9b) is false.  These results include 
the previous results on primes and sums of two squares, and also cover 
many other examples.  

Consider sets containing roughly half of the prime numbers.  There are uncountably many 
such sets, and so maybe we can find a set which is very well distributed in arithmetic 
progressions.  One amusing example from [16] shows that this cannot be done, and the 
Friedlander-Granville limitations persist for any such set.  

We content ourselves with this vague description 
of the uncertainty principle,  referring the reader to [16] for more examples 
and a precise description of the results.


\Refs

 \ref\key 1 
 \by A. Balog and T.D. Wooley 
 \paper Sums of two squares in short intervals
 \jour Canad. J. Math. 
 \vol 52
 \yr 2000
 \pages 673--694
 \endref

 \ref\key 2
 \by E.B. Bogomolny and J.P. Keating
 \paper Random matrix theory and the Riemann zeros. II. $n$-point correlations
 \jour Nonlinearity
 \vol 9
 \pages 911--935
 \yr 1996
 \endref

 \ref\key 3
 \by T.H. Chan 
 \paper Pair correlation and distribution of prime numbers 
 \jour Ph. D. Thesis, University of Michigan 
 \yr 2002
 \pages 101pp
 \endref

 \ref\key 4 
 \by T.H. Chan 
 \paper A note on primes in short intervals
 \jour Int. J. Number Theory
 \vol 2 
 \yr 2006
 \pages 105--110
 \endref

 \ref\key 5
 \by H. Cram{\' e}r 
 \paper On the order of magnitude of the difference between consecutive prime numbers 
 \jour Acta Arith.
 \vol 2
 \pages 23--46
 \yr 1936
 \endref

 \ref\key 6
 \by H. Davenport 
 \book Multiplicative number theory
 \publ Springer Graduate Texts in Math. 74 
 \yr 2000
 \endref

 \ref\key 7
 \by W. Feller
 \book An introduction to probability theory and its applications
 \publ Wiley
 \yr 1966 
 \endref

 \ref\key 8
 \by J. Friedlander and A. Granville 
 \paper Limitations to the equi-distribution of primes I 
 \jour Annals of Math. 
 \vol 129
 \yr 1989 
 \pages 363--382
 \endref

\ref\key 9
\by P. X. Gallagher 
\paper On the distribution of primes in short intervals 
\jour Mathematika 
\vol 23
\yr 1976
\pages 4--9
\endref

\ref\key 10
\by D. Goldston 
\paper Notes on pair correlation of zeros and prime numbers
\inbook Recent perspectives in random matrix theory and number theory 
\publ London Math. Soc. Lecture Notes Ser. 322, Cambridge U. Press 
\yr 2005
\pages 79--110
\endref

\ref\key 11
\by D. Goldston and H.L. Montgomery 
\paper On pair correlations of zeros and primes in short intervals 
\inbook Analytic Number Theory and Diophantine Problems 
\publ Prog. in Math. Birkh{\"a}user \vol 70
\yr 1987
\pages 183--203
\endref

\ref \key 12 
\by D. Goldston, J. Pintz and C. \Yildirim 
\paper Primes in tuples, I 
\jour Ann. of Math. (to appear), preprint available at www.arxiv.org
\endref 


\ref\key 13
\by A. Granville 
\paper Unexpected irregularities in the distribution of prime numbers
\jour Proc. of the Int. Congr. of Math., Vol. 1, 2 (Z{\"u}rich, 1994) 
\pages 388-399
\publ Birkh{\" a}user, Basel
\yr 1995   
\endref

\ref\key 14
\by A. Granville and G. Martin 
\paper Prime number races
\jour Amer. Math. Monthly
\vol 113
\yr 2006
\pages 1--33
\endref

\ref\key 15
\by A. Granville and K. Soundararajan 
\paper Sieving and the Erd{\H o}s-Kac theorem
\inbook these proceedings
\endref

\ref\key 16
\by A. Granville and K. Soundararajan 
\paper An uncertainty principle for arithmetic sequences 
\jour Ann. of Math. (to appear), preprint available at www.arxiv.org
\endref 

\ref\key 17
\by G.H. Hardy and J.E. Littlewood 
\paper Some problems of Paritio Numerorum (III): On the expression of a number as a 
sum of primes 
\jour Acta Math. 
\vol 44
\pages 1--70
\yr 1922
\endref


\ref\key 18
\by D.R. Heath-Brown
\paper Differences between consecutive primes 
\jour Jahresber. Deutsch. Math.-Verein. 
\vol 90
\yr 1998
\pages 71--89
\endref

\ref\key 19
\by A. Hildebrand and H. Maier
\paper Irregularities in the distribution of primes in short intervals
\jour J. Reine Angew. Math. 
\vol 397
\yr 1989 
\pages 162--193
\endref

\ref\key 20
\by C. Hooley
\paper On the difference between consecutive numbers prime to $n$: II
\jour Publ. Math. Debrecen 
\vol 12
\yr 1965 
\pages 39--49
\endref

\ref\key 21
\by C.P. Hughes and Z. Rudnick 
\paper On the distribution of lattice points in thin annuli
\jour Int. Math. Res. Not.
\vol 13 \yr 2004 
\pages 637--658
\endref

\ref\key 22
\by H. Maier
\paper Primes in short intervals
\jour Michigan Math. J. 
\vol 32
\yr 1985
\pages 221--225
\endref

\ref\key 23
\by W.R. Monach 
\paper Numerical investigation of several problems in number theory
\jour Ph. D. Thesis, University of Michigan
\yr 1980 
\pages 171 pp 
\endref


\ref\key 24
\by H.L. Montgomery 
\paper The pair corelation of zeros of the zeta function 
\inbook Analytic Number Theory (St. Louis Univ. 1972) 
\publ Proc. Sympos. Pure Math. (Amer. Math. Soc.) 
\vol 24
\pages 181-193
\yr 1973
\endref

\ref\key 25 
\by H.L. Montgomery and K. Soundararajan
\paper Beyond pair correlation 
\inbook Paul Erd{\H o}s and his mathematics, I 
\publ Bolyai Soc. Math. Stud., 11, Budapest
\yr 2002
\pages 507-514 
\endref


\ref\key 26
\by H.L. Montgomery and K. Soundararajan 
\paper Primes in short intervals 
\jour Comm. Math. Phys. 
\vol 252
\yr 2004 
\pages 589--617
\endref

\ref\key 27
\by H.L. Montgomery and R.C. Vaughan
\paper On the distribution of reduced residues 
\jour Annals of Math.
\vol 123
\yr 1986
\pages 311--333
\endref

\ref \key 28
\by E. M. Rains
\paper High powers of random elements of compact Lie groups
\jour Probab. Theory Related Fields \vol 107 \yr 1997 \pages 219--241
\endref



\ref\key 29
\by M. Rubinstein and P. Sarnak 
\paper Chebyshev's Bias
\jour Experimental Math.
\vol 3
\yr 1994 
\pages 173--197
\endref

\ref\key 30
\by A. Selberg
\paper On the normal density of primes in short intervals, and the difference 
between consecutive primes
\inbook Collected papers (Volume I) 
\publ Springer
\yr 1989 
\pages 160--178
\endref

\ref \key 31 
\by K. Soundararajan 
\paper Small gaps between prime numbers: the work of Goldston-Pintz-\Yildirim
\jour Bull. Amer. Math. Soc. (to appear), preprint available at www.arxiv.org
\endref


 \endRefs
 \enddocument